\documentclass[12pt]{article}
\usepackage{latexsym}

\newtheorem{theorem}{Theorem}[section]
\newtheorem{cor}[theorem]{Corollary}
\newtheorem{lemma}[theorem]{Lemma}
\newtheorem{remark}[theorem]{Remark}
\newtheorem{defin}[theorem]{Definition}

\newtheorem{prop}[theorem]{Proposition}

\def\reff#1{{\rm
(\ref{#1})}}

\def\R{{\rm I\kern-.25em R}}
\def\F{{\rm I\kern-.25em F}}
\def\C{{\rm \kern.24em \vrule width.02em height1.4ex
depth-.05ex \kern-.26em C}}

\newcommand{\ms}{\nu^{0}_{H}} 
\newcommand{\leqs}{ \ \leq \ }

\newcommand{\CC}{{\bf (C)}}
\newcommand{\DD}{{\bf (D)}}

\newcommand{\be}{\begin{equation}}
\newcommand{\ee}{\end{equation}}
\newcommand{\bd}{\begin{displaymath}}
\newcommand{\ed}{\end{displaymath}}
\newcommand{\bp}{\underline{\bf Proof}:\ }
\newcommand{\ep}{{\hfill $\Box$}\\}
\newcommand{\ba}{\begin{array}}
\newcommand{\ea}{\end{array}}
\newcommand{\e}{ \ = \ }

\newcommand{\w}{ \omega }

\newcommand{\gla}{\gamma_{\lambda}}
\newcommand{\ala}{\alpha_{\lambda}}

\begin{document}
\title{PARTIAL NORMS AND THE CONVERGENCE OF
GENERAL PRODUCTS OF MATRICES}
\author{Michael Neumann\thanks{Work supported by NSF Grant
DMS-9306357.}
\\Department of Mathematics
\\University of Connecticut
\\ Storrs, CT 06268, USA
\and Hans Schneider
\footnote {Research supported in part by NSF Grant DMS-9424346.
This research was partly performed during a visit by this author
to the University of Bielefeld, Germany, supported by the
Sonderforschungsbereich ``Diskrete Strukturen in der
Mathematik'', Universit\"at Bielefeld.}
\\Department of Mathematics
\\University of Wisconsin
\\ Madison, WI 53706, USA}

\date{21 November 1997}

\maketitle

\begin{abstract}
Motivated by the theory of inhomogeneous Markov chains, we
determine
a sufficient condition for the convergence to $0$ of a general
product
formed from a sequence of
real or complex matrices. When the matrices have a
common invariant subspace $H$, we give a sufficient condition
for the convergence to $0$ on $H$ of a general product. Our
result
is applied to obtain a condition for the weak ergodicity of an
inhomogeneous Markov chain.
We compare various types of contractions
which may be defined for a single matrix, such as
paracontraction,
$l$--contraction, and $H$--contraction,
where $H$ is an invariant subspace of the matrix.

\end{abstract}

\thispagestyle{empty}

\newpage

\section{Introduction}
\label{sec1}
\setcounter{equation}{0}

Recently there has been much interest in conditions for the
convergence
of infinite products of real or complex matrices. Several
investigations have concentrated
on products taken in one direction -- left or right,
see for example the recent papers by Beyn and Elsner
\cite{BeEl97} and Hartfiel and Rothblum \cite{HaRo97}.
However, in this paper,
we are concerned with {\it general products} formed from a
given infinite sequence
of matrices. These are defined further on in the paper
and they have previously been considered for nonnegative and for
stochastic matrices by Seneta in \cite[Chapters 3.1 and
4.6]{Sene81}.\\

Our principal result is a sufficient condition for the
convergence to $0$ of infinite general products of matrices.
We pay particular attention to the case where
there is a common invariant subspace for all the matrices
in the product. As a special case, we obtain a result on the
weak ergodicity of an inhomogeneous Markov chain.\\

The investigations described above  are preceded
by a study of the
interrelations of several types of contractions which may be
defined for a single matrix.\\

The  motivation for our study comes from the theory of
inhomogeneous
Markov chains, see \cite {Sene73} and \cite[Chapter 4]{Sene81}
for
much further
background material. This theory naturally leads to the study
of the convergence, ergodicity and weak ergodicity of products of
stochastic matrices.\\

We now describe our paper in more detail.
First, in Section 2, we  sharpen several
observations
concerning paracontraction. This we shall achieve
by introducing a partial norm on the matrix restricted to an
invariant
subspace $H$ and the consequent notion
of $H$--{\it contraction}. We relate this concept to
the notion of paracontraction as introduced by Nelson and
Neumann in
\cite{NeNe87}  and to the notion of $l$--contraction
recently introduced by Beyn and Elsner in  \cite{BeEl97}.\\

Second, in Section 3,
which is essentially independent of Section 2,
we turn to our main results.
We formulate our sufficient conditions depending on norms
for the convergence to $0$ of infinite general products of
matrices
which are formed from a given
infinite sequence of matrices.
We then apply the results to examine the convergence to $0$
of infinite general
products on a common invariant subspace of the matrices.\\

In Section 4 we specialize
our results and we investigate
products of stochastic matrices to obtain a suffcient condition
for weak
ergodicity. The $\ell_1$ coefficient
of ergodicity due to Bauer, Deutsch, Stoer \cite{BaDeSt69} plays
a
special role here.\\

In Section 5
we give bounds on the algebraic and geometric multiplicities of
the eigenvalues of the restriction of a matrix to an invariant
subspace in terms of these quantities for the entire matrix.
These are related to familiar results on stochastic matrices.

\section{Partial norms and Paracontractions}
\label{sec2}

In this paper $\F$ will stand for the real field $\R$ or the
complex
field  $\C$.
We begin by recalling  the following definition:\\

\begin{defin}
\label{parcont} {\rm (Nelson and Neumann \cite{NeNe87}) Let $A\in
\F^{n,n}$
and let $\nu$
be  a vector norm on $\F^{n}$. Then $A$ is {\it paracontracting
with
respect to} $\nu$ if
\be
\label{para.cond}
  \nu(Ax) \ < \nu(x) \ \ {\rm whenever } \ Ax \neq x.
\ee
}
\end{defin}

In particular we see that one immediate consequence of the
definition is that
$\nu^0(A)\leq 1$,
where $\nu^0(A)$ denotes the operator norm corresponding to the
norm
$\nu$.
For a given matrix norm let ${\cal N}_{\nu}$ denote {\rm
the set of
all paracontracting matrices with respect to} $\nu$ in $\F^{n}$.
Nelson and Neumann show in \cite{NeNe87} that if $A\in {\cal
N}_{\nu}$, then
$\lim_{i\rightarrow
\infty}A^i$ exists and that
$AB\in {\cal
N}_{\nu}$
if $A,B\in {\cal N}_{\nu}$.
\\

We now introduce the concept of a partial norm with respect to a
subspace
and the concept of H--contraction:\\

\begin{defin}
\label{Hcont}
{\rm
Let $H$ be a subspace of $\F^n$ and let $A \in \F^{nn}$ be
a matrix which leaves $H$ invariant.
Let $\nu$ be a norm on $\F^n$ and
define $\nu^0$ to be the operator norm induced by $\nu$.
Then
\be
\label{Hcond}
\nu^{0}_{H}(A) \ = \ \sup_{0 \neq x\in
H}\frac{\nu(Ax)}{\nu(x)}.
\ee
We call $\nu^{0}_{H}(A)$ a {\em partial norm of} $A$ {\em with
respect
to the invariant subspace} $H$ of $A$. }
\end{defin}

\begin{defin}
\label{nonexpan}
{\rm
 The matrix  $A\in \F^{n,n}$ is {\em nonexpansive} with respect
 to the
vector norm
$\nu$ on $\F^{n}$ if
$\nu^{0}(A)\leq 1$. The matrix  $A$
is an  $H$--{\em
contractor} if
it is nonexpansive, $H$ is an invariant subspace of $A$, and
$\nu_{H}^{0}(A)<1$. }
\end{defin}

We shall
denote the range of $A \in \F^{nn}$ by ${\cal R}(A)$ and the
nullspace
of $A$ by ${\cal N}(A)$.

\begin{remark}
\label{complement}
{\rm If $A\in \F^{n,n}$ is nonexpansive with respect to $\nu$
then
the spectral  radius $\rho(A)$ of $A$ satisfies
$\rho(A)\leq \nu^0(A) \leq 1$. If $1$ is an
eigenvalue of $A$, then we must have that $\nu^0(A)=\rho(A)$ and
it follows that all Jordan blocks of $A$
corresponding to $1$ are
$1\times 1$, (viz. index$_0(A) = 1$), e.g., Mott and Schneider
\cite{MoSc59}. Hence $K={\cal{N}}(I-A)$ and
$H={\cal R}(I-A)$ are
complementary invariant subspaces.}
\end{remark}

We next show the connection between paracontraction and
$H$--contraction.

\begin{lemma}
\label{lemma1}
Let $A\in \F^{n,n}$ be paracontracting with respect to
$\nu$ and let $K$ be the subspace of all its fixed points.
Then $H={\cal R}(I-A)$ is
invariant under $A$, complementary to $K$, and such that
$\nu^{0}_{H}(A)<1$
so that
$A$ is an
$H$--contractor.\\

\end{lemma}

\bp
It follows by Remark \ref{complement} that $H$ is invariant under
$A$
and complementary to $K$.
Since $A$ is paracontracting with respect
to $\nu$, we have that  $\nu^0(A) \leq 1$.
In view
of \reff{para.cond}, it now follows that
$\nu(Ax)<\nu(x)$, for all
$ x \in H$, and so
since the
unit ball
${\cal U}$ of $\nu^{0}_{H}$ in $H$ is compact,
\bd
\max_{x\in {\cal U}}\nu(Ax) \ < \ 1.
\ed
Hence $\nu_{H}^{0}(A)<1$.
\ep

Recently Beyn and Elsner \cite{BeEl97} have introduced the notion
of
$l$--paracontraction:\\

\begin{defin}
\label{elsner}
{\rm (\cite{BeEl97}) A matrix $A \in \F^{n,n}$ is
$l$--paracontracting
with respect to the vector norm $\nu$ on $\F^{n}$ if there exists
a
$\gamma>0$ such that
\be
\label{beyn.elsner}
   \nu(Ax) \leqs \nu(x)-\gamma \nu(Ax-x), \ \ \forall x\in
   \F^{n}.
\ee
}
\end{defin}

In their paper Beyn and Elsner establish
several conditions which are equivalent to $l$--paracontraction.
Clearly $l$--paracontraction implies paracontraction.
We next
show that
H--contraction implies $l$--paracontraction
for a suitable chosen norm and a suitably chosen subspace $H$.\\

\begin{theorem}
\label{equiv}
Let $A\in \F^{n,n}$ and suppose that $A$ is nonexpansive
with respect to the norm
$\nu$
on $\F^n$.
Let $H={\cal R}(I-A)$ and let $K = {\cal N}(I-A)$.
If the norm $\nu$ satisfies
\begin{equation}
\label{additive}
\nu(y+z) = \nu(y) + \nu(z), \ {\rm for} \ y \in H, \ z \in K,
\end{equation}
then the following are equivalent:\\

{\rm (i)} $A$ is $l$--paracontracting with respect to $\nu$.\\

{\rm (ii)} $A$ is paracontracting with respect to $\nu$.\\

{\rm (iii)} $A$ is an $H$--contraction with respect to $\nu$.
\end{theorem}

\bp
(i) $\Longrightarrow$ (ii): Obvious.\\

(ii) $\Longrightarrow$ (iii): By Lemma \ref{lemma1}.\\

(iii) $\Longrightarrow$ (i): Since $A$ is nonexpansive, we have
by
Remark \ref{nonexpan} that $K \oplus H \ = \F^n$.
For $x \in \F^n$, write
\bd
x \ = \ y + z, \ \ \ \ y \in H, \ z \in K.
\ed
Let $k = \nu^0_H(A)$ and note that $k < 1$ by (iii).
Since $x - Ax$ = $y - Ay \in H$ it follows that
\be
\label{inequ1}
\nu(x-Ax) = \nu((I-A)y) \ \leq \ \nu^0_H(I-A)\nu(y) \
\leq \ (1+k)\nu(y).
\ee
Moreover, using \reff{additive}, we obtain
\be
\label{inequ2}
\nu(x) -
\nu(Ax) \ = \
\nu(y+z)-
\nu(Ay+z) \ = \
\nu(y)-\nu(Ay) \ \geq \ (1-k)\nu(y).
\ee
We combine \reff{inequ1} and \reff{inequ2} and we deduce that
\begin{displaymath}
\nu(x-Ax) \ \leq \ (1+k)\nu(y) \ \leq \ \frac{1+k}{1-k}
(\nu(x)-\nu(Ax)).
\end{displaymath}
and this is equivalent to \reff{beyn.elsner} with
$ \gamma = (1-k)/(1+k)$.
\ep

\begin{remark} {\rm
If $\nu'$ is a norm on $\F^n$, then we may define a norm $\nu$
that satisfies \reff{additive} and agrees with $\nu'$ on $H$ and
$K$
by setting
$\nu(y+z) = \nu'(y) + \nu'(z)$  for $y \in H$ and $z \in K$.}
\end{remark}
\begin{remark} {\rm
The implication (iii) $\Longrightarrow$ (ii) of Theorem
\ref{equiv}
holds under weaker assumptions on the norm $\nu$. Specifically,
if
$A$ is an $H$--contraction with respect to $\nu$ and the
condition
\bd
\nu(y) < \nu(y') \Rightarrow \nu(y+z) < \nu(y'+z), \ y \in H,\ z
\in K
\ed
is satisfied, then $A$ can be shown to be paracontracting
with respect to $\nu$.}
\end{remark}

\section{Convergence of Infinite Products}
\label{sec3}

In this section we develop our main results concerning the
 convergence of  products of  complex
matrices taken in an arbitrary order from an infinite sequence of
matrices.
Such products were considered (in a slightly less general form)
in Seneta \cite[Section 4.6]{Sene81}
in the case of stochastic matrices, see also \cite{Leiz92} and
\cite{Rhod97}.\\

Let $A_1,A_2, \ldots$
be a sequence of complex matrices. We shall consider
products of matrices  obtained from the sequence in the following
manner: First choose
some permutation of the given infinite sequence to obtain a
sequence
$B_1,B_2,\ldots$. Then form the products $C_{p,r}$
 of the matrices $B_{p+1},
\ldots,
B_{p+r}$ in some order. We shall call $C_{p,r}$ a {\em general
product} from
the sequence $A_1,A_2,\ldots$
and we shall consider the existence of
$\lim_{r\rightarrow\infty} C_{p,r}$.
If this limit is $0$,
for all permutations of $A_1, A_2, \ldots $ and all $p, \ p \geq
0$,
then we shall say that all general products from the
the sequence $A_1,A_2,\ldots$ converge to $0$.\\

As an example of a sequence of general product suppose the chosen
order is
$A_9, A_7, A_5$, $A_{14}, A_2, \ldots$. Then
the sequence of $(C_{2,1}, C_{2,2}, \ldots)$ may
begin thus:
\bd
  \begin{array}{l}
          C_{2,1} \e A_7, \ \\
                          \ \\
          C_{2,2} \e A_7 A_5, \ \\
                      \ \\
          C_{2,3} \e A_5 A_2 A_7 ,\ \\
                                       \ \\
          C_{2,4} \e A_2 A_7 A_{14} A_5.
  \end{array}
\ed
Note that, for a given sequence $C_{p,1}, C_{p,2}, \ldots$
of general products each factor of $C_{p,r}$ occurs in
$C_{p,r+1}$,
but the order in which the factors occur in $C_{p,r}$ is
arbitrary.\\

Let $\mu$ be a matrix norm (viz. a submultiplicative norm
on $\F^{nn}$)
and denote
\bd
   \mu^+(P) \e \max (\mu(P), 1) \;\;\;\;\;\;
   \mbox{and}\;\;\;\;\;\;
   \mu^-(P) \e \min (\mu(P), 1).
\ed

Now let $A_1,A_2,\ldots$ be a sequence of
matrices in $\F^{nn}$ and let $\mu$ be a  matrix norm.
We now define two conditions:\\
\begin{description}
\item[{\fbox{Condition (C)}}]: We say that
the sequence $A_1,A_2,\ldots,$
satisfies Condition\\

 {\bf (C)}  for the norm $\mu$ if
\be
\label{condC}
\sum_{i=1}^{\infty} (\mu^+(A_i)-1) \ {\rm converges}.
\ee

\item[{\fbox{Condition (D)}}]: We say that the sequence
$A_1,A_2,\ldots$
satisfies
Condition\\

 {\bf (D)} for the norm $\mu$ if
\be
\label{condD}
   \sum_{j=1}^{\infty} (1-\mu^-(A_i)) \
             {\rm diverges. }
\ee
\end{description}

We are now ready to prove the following result:\\

\begin{prop}
\label{prop.bdd}
Let $A_1,A_2,\ldots$ be a sequence of  matrices in $\F^{nn}$.
Let $\mu$ be a matrix norm on $\F^{nn}$.  Suppose
that the sequence $A_1,A_2,\ldots$
satisfies Condition \CC \ for the norm $\mu$.
 Then all general products from $A_1,A_2,\ldots$ are bounded.
 \end{prop}

\bp
Let $B_1,B_2, \ldots$ be a permutation of $A_{1}, A_2, \ldots $
and
let $C_{p,r}$ be a product of $B_{p+1}, \ldots, B_{p+r}$ in some
order.
By  Condition \ \CC  and   \cite[Theorem 14]{Hysl45},
$\sum_{i=1}^{\infty}(\mu^+(B_i)-1)$
converges and hence
$\sum_{i=1}^{\infty}(\mu^+(B_{p+i})-1)$ also converges.
Thus, by
\cite[Theorem 51]{Hysl45}, the product
$\prod_{i=1}^{\infty} \mu^+(B_{p+i})$
converges and so there exists a positive constant
$M$ such that $\prod_{i=1}^r \mu^+(B_{p+i})\leq M$, for each $r
\in
\{1,2,\ldots\}$.
It follows that
\begin{equation}
\label{eqn1}
   \begin{array}{lll}
 \mu\left(C_{p,r}\right ) & \leq & \mu\left( B_{p+1} \right)
 \cdots
\mu\left(B_{p+r}\right)
\\
                   & \ & \ \\
                   &  =   & \left[\prod_{i=1}^{r}
\mu^-\left(B_{p+i}\right)\right]
                            \left[\prod_{i=1}^{r}
\mu^+\left(B_{p+i}\right)\right]
\\
                    \ & \ & \ \\
                   & \leq & M \left[\prod_{i=1}^{r}
\mu^-\left(B_{p+i}\right)\right]
\\
                  \ & \ & \ \\
				   & \ \leq  & M.
   \end{array}
\end{equation}

\ep

The above proposition allows us to prove a stronger result under
an additional conditions. Note that in the theory of infinite
products
of nonnegative numbers it is customary to speak of {\it
divergence} to
$0$, see e.g. \cite[p. 93]{Hysl45}. \\

\begin{theorem}
\label{theorem.conv}
Let $A_1,A_2,\ldots$ be a sequence of  matrices in $\F^{nn}$.
Let $\mu$ be a matrix norm on $\F^{nn}$.  Suppose
that the sequence $A_1,A_2,\ldots$
satisfies Conditions \CC \ and \DD \
for the norm $\mu$.
 Then all general products from $A_1,A_2,\ldots$ converge to $0$.
 \end{theorem}

\bp
Let $B_1,B_2, \ldots$ be a permutation of $A_{1}, A_2, \ldots $
and
let $C_{p,r}$ be a product of $B_{p+1}, \ldots, B_{p+r}$ in some
order.
As in the proof of Proposition(\ref{prop.bdd}), we have that
\begin{equation}
\label{eqn2}
 \mu\left(C_{p,r}\right )
                  \  \leq \
				   M \left[\prod_{i=1}^{r}
\mu^-\left(B_{p+i}\right)\right].
\end{equation}

By Condition \DD,
the sum  $\sum_{i=1}^{\infty} (1-\mu^-(A_i))$
diverges and so, by \cite[Theorem 14]{Hysl45},
$\sum_{i=1}^{\infty} (1-\mu^-(B_i))$ diverges.
Thus
$\sum_{i=1}^{\infty} (1-\mu^-(B_{p+i}))$ also diverges.

We again apply  \cite[Theorem 51]{Hysl45}, and we obtain that
$\prod_{i=1}^{\infty} \mu^-(B_{p+i})$ diverges.
But since $\mu^-(B_i) \leq 1$, the last product must diverge to
$0$
and the proof is done.
\ep

If $H$ is a subspace of $\F^n$ which is invariant under  $A \in
\F^{nn}$,
we denote the restriction of $A$ to $H$ by $A_{|_H}$.
As an immediate application of Theorem \ref{theorem.conv} we
obtain the
following result:\\

\begin{cor}
\label{ergod.prod}
 Let $A_1,A_2,\ldots$ be a sequence of matrices
 and let  $C_{p,1},
C_{p,2},\ldots$ be
sequence of general products from the $A_1,A_2,\ldots$.
Let $H$ be a subspace of $\F^n$ which is invariant
under each $A_i, \  i = 1,\ldots$. Let $\nu$ be a norm on $\F^n$.
Suppose that  the sequence $((A_i)_{|H}) $ satisfies conditions
\CC  \
and \DD\  for the norm $\nu^0_H$.
If $x \in H$, then
\bd
\lim_{r \rightarrow \infty} C_{p,r}x  \e 0 .
\ed
\end{cor}

\bp
 Immediate by Theorem \ref{theorem.conv}.
\ep

Since $\nu^0_H(A) \leq \nu^0(A)$ when $H$ is an invariant
subspace of
$A$, it follows easily that under the conditions of
Corollary(\ref{ergod.prod}),
Condition \CC (resp. Condition \DD)
on the sequence of $(A_i)$ implies
Condition \CC (resp. Condition \DD)
on the sequence of
$((A_i)_{|H})$.

\section{Applications to Stochastic Matrices}
\label{sec4}

In this section we  apply the foregoing results
to stochastic matrices. In
order to be consistent with our previous section
we consider {\em column} stochastic matrices.Thus ``stochastic
matrix''
will mean ``column stochastic matrix''.\\

Let $e = (1,\ldots,1)^T \in \R^n$. Hence
 in this section we shall assume
that
\be
\label{partic.H}
H \e \{x \in \R^n \ : \ e^Tx= 0 \}.
\ee
If $A$ is a stochastic matrix in $\R^{n,n}$, then $H$
is invariant under $A$, but note that $H$ need not be a
complement
to the space of fixed points of $A$.
If $\nu$ is a norm on $\R^n$ and $H$ is given by \reff{partic.H},
then
the corresponding partial norm $\nu_H^0$ on
$\R^{nn}$ will be called a {\em coefficient of ergodicity} as is
usual in the
literature on Markov chains.\\

The $\ell_1$ norm on $\R^n$ plays a special
role in this theory and we shall denote it henceforth by $\w$.
The corresponding coefficient of ergodicity was apparently first
computed
in Bauer, Deutsch, and Stoer \cite{BaDeSt69}, see also
\cite{Zeng72},
and equals
\bd
\w_{H}^{0}(A) \e (1/2) \max_{i,k \in \{1,\dots,n\}} \sum_{j=1}^n
|a_{i,j}-a_{k,j}|.
\ed
It is known that $\w_{H}^{0}(A) \leq 1$ for all stochastic
matrices $A$
and $\w_{H}^{0}$
is the only coefficient of ergodicity that satisfies this
inequality,
\cite{KuRh90} and \cite{Lesa90},
but see also \cite{Rhod97}.\\

\begin{remark} {\rm
We comment that the only
stochastic matrix $A$ which is paracontracting with respect to
the $\ell_1$
norm $\w$
is the
identity matrix. This follows easily from the facts
that for any nonnegative vector $u$ with
$\w(u) = 1$, we also have that $\w(Au) = 1$.
On the other hand, every stochastic matrix that has a
row with two or more
positive elements
is H--contracting
for the subspace $H$ given
in \reff{partic.H} and the norm $\w$.
}\end{remark}

\begin{defin}
\label{weak}
{\rm
Let $P_1,P_2,\ldots$ be a sequence of $n \times n$ stochastic
matrices.
We shall say
that all general products
formed from
this sequence are {\em weakly ergodic} if for all general
products
$B_{p,1}, B_{p,2}, \ldots$, we have that
\be
\label{weak.con}
lim_{r \rightarrow \infty}B_{p,r}x \e 0, \ \  {\rm for \ all} \ \
x \in H.
\ee.
}
\end{defin}

Since every $x \in H$ can be written $x = c(u - v)$, where $u,v
\in \R^n$
are nonnegative and $e^Tu = e^Tv = 1$ and $c \in \R$, it is
easily seen
that, for each product considered, our definition is equivalent
to that in
\cite{Hajn58}, \cite{MoSc57}, or \cite[Defn. 3.3]{Sene81}.\\

By Theorem \ref{theorem.conv} we now immediately obtain:\\

\begin{theorem}
\label{wkerg}
Let $\nu$ be a norm on $\R^n$ and let
$\ms$ be the corresponding coefficient of ergodicity.
Let $P_1,P_2,\ldots$ be a sequence of $n \times n$ stochastic
matrices.
Then all general products formed from
this sequence  are weakly ergodic  if

\be
\label{condCW}
\sum_{i=1}^{\infty} ({\ms}^+(P_i)-1) \ {\rm converges}
\ee
and
\be
\label{condDW}
   \sum_{i=1}^{\infty} (1-{\ms}^-(P_{i_j})) \
             {\rm diverges \ for\  some \  subsequence\ of} \
             P_1,P_2\ldots.
\ee
\end{theorem}

Note that (\ref{condCW}) is automatically satisfied if $\nu =
\w$,
the $\ell_1$-norm. This special case of Theorem \ref{wkerg} is
observed in \cite[Exercise 4.36]{Sene81}.
Results related to this case (and which therefore
involve only (\ref{condDW}) explicitly) are to be found in
\cite[Theorem 1]{Sene73}  and in \cite{MoSc57}. The theorem in
the
latter paper is
there illustrated by an
example of a sequence of stochastic matrices
that satisfies (\ref{condDW}) for
the norm $\w$, see \cite[p. 333]{MoSc57}.\\

The following corollary is due to  Rhodius
\cite[Thm. 3, Part I]{Rhod97} in the case of $\nu = \w$,
see Leizarowitz \cite[Thm A (i)]{Leiz92} for a
related result.

\begin{cor}
\label{accum}
Let $P_1, P_2, \ldots$ be a sequence of stochastic matrices. If
$\ \nu_H^0(P_i)$ \  $\leq 1$ for all $i, \ i = 1,2, \ldots$, and
there exists a point of accumulation $c$ of the sequence
$\nu_H^0(P_1),\nu_H^0(P_2), \ldots$ such that $c < 1$, then all
general
products of the sequence are weakly ergodic.
\end{cor}

\bp
Clearly condition \reff{condCW} holds and
there is an infinite
subsequence $j_1, j_2, \ldots$ such that $\w_H^0(P_{i_j}) <
(1+c)/2 < 1$,
$j = 1,2, \ldots$. Then condition \reff{condDW} holds for this
subsequence.
The result follows from Theorem (\ref{wkerg}).
\ep

Another corollary of Theorem \ref{wkerg} is \cite[Thm A
(ii)]{Leiz92}:

\begin{cor}
\label{accum2}
Let $P_1, P_2, \ldots$ be a sequence of stochastic matrices and
let
$\nu$ be a norm on $\R^n$.
If all points of accumulation $c$ of the sequence
$\nu_H^0(P_1)$, \ $\nu_H^0(P_2), \ldots$ satisfy $c < 1$, then
all general
products of the sequence are weakly ergodic.
\end{cor}

\bp
Since the set of accumulation points of a bounded sequence is
compact,
there exists $d < 1$ such that only a finite number of terms of
the
sequence $\nu_H^0(P_1),\nu_H^0(P_2), \ldots$ exceed $d$. Hence
\reff{condCW}
and \reff{condDW} hold for the sequence of ergodicity
coefficients
and the corollary follows from Theorem \ref{wkerg}.
\ep

\section{Bounds for eigenvalues}
\label{sec5}

\begin{theorem}
\label{thm.bound}
Let $H$ be an invariant subspace for $A \in \C^{n,n}$. Let
$\ala(A)$
and $\gla(A)$ be the algebraic and geometric
multiplicities of $\lambda \in \C^{nn}$ as an eigenvalue of  $A$,
respectively. Then
\bd
\ala(A) \ \leq \ \ala(A_{|_H}) + n - \dim(H)
\ed
and
\bd
\gla(A) \ \leq \ \gla(A_{|_H}) + n - \dim(H).
\ed

Further, if $\mu$ is a matrix norm and
$|\lambda| > \mu_H^0(A)$, then $\ala(A) \leq n - \dim(H)$ and
$\gla(A) \leq n - \dim(H)$.
\end{theorem}

\bp
By a proper choice of basis for $\C^n$ we may put $A$ in the form
\be
\label{eqn3}
  A \e \pmatrix{A_{1,1} & 0 \cr A_{2,1} & A_{2,2} \cr}
\ee
and the first inequality follows immediately.\\

The second inequality is an immediate consequence of
\cite[Theorem
2.1]{MeRo77}, see also \cite[Cor. 5.8]{HeRoSc89}.\\

Finally, since all eigenvalues $\lambda$ of $A_{22}$ satisfy
$|\lambda| \leq \mu_H^0(A)$, the last part of the theorem
follows.
\ep

\begin{remark}
{\rm
Theorem (\ref{thm.bound}) implies several known results.
We now give an example below which could also be derived from
\cite[Theorem 2]{HaRo97}. Let $A$ be a stochastic matrix, let $H$
be
defined by
\reff{partic.H},
and suppose that a coefficient of ergodicity
satisfies $\nu^{0}_{H}(A) < 1$. Since, in this case, the matrix
$A_{1,1}$ in (\ref{eqn3}) is $1 \times 1$, it follows that
$1$ is an
algebraically simple eigenvalue of $A$ and all other eigenvalues
$\lambda$ of $A$ satisfy  $\nu^{0}_{H}(A) \geq |\lambda|$. }
\end{remark}

\centerline{{\bf Acknowledgment}}
\vspace{.15in}

We would like to thank Wenchao Huang for some helpful remarks.
We thank Olga Holtz for her careful reading of the manuscript.\\

\bibliographystyle{plainxe}

\begin{thebibliography}{0}


\bibitem{BaDeSt69}
F. L. Bauer, Eck. Deutsch, and J. Stoer.
\newblock Absch\"{a}tzungen f\"{u}r die Eigenwerte
positiver linearer Operatoren.
\newblock {\it Lin. Alg. Appl.}, 2:275--301, 1969.

\bibitem{BeEl97}
W. J. Beyn and L. Elsner.
\newblock Infinite products and paracontracting matrices.
\newblock {\it Elec. J. Lin. Alg.}, 2:1--8, 1997.

\bibitem{Hajn58}
J. Hajnal,
\newblock Weak ergodicity in non-homogeneous Markov chains,
\newblock {\it Proc. Cambridge Phil. Soc.}, 54:52--67, 1956

\bibitem{HaRo97}
D.J. Hartfiel and U.G. Rothblum,
\newblock Convergence of inhomogeneous products of matrices and
coefficients of ergodicity.
\newblock {\it Lin. Alg. Appl.}, to appear.

\bibitem{HeRoSc89}
D. Hershkowitz, U. G. Rothblum, and H. Schneider.
\newblock The combinatorial structure of the generalized
nullspace
of a block triangular matrix.
\newblock {\it Lin. Alg. Appl.}, 116:9--26, 1989.

\bibitem{Hysl45}
J. M. Hyslop.
\newblock Infinite Series.
\newblock {\it Oliver and Boyd}, Edinburgh, 1945.

\bibitem{KuRh90}
R. K\"uhne and A. Rhodius, A characterization of
Dobrushin's coefficient of ergodicity,
\newblock {\it Z. Anal. Anw.} 9(2):187--188,
1990.

\bibitem{Leiz92}
A.Leizarowitz.
\newblock On infinite products of stochastic matrices,
\newblock {\it Lin. Alg. Appl.} 168:189--219, 1992.

\bibitem{Lesa90}
A. Le\v{s}anovsk\'{y},
\newblock Coefficients of ergodicity generated by non-symmetrical
vector norms,
\newblock {\it Czechoslovak Math. J.} 40 (115):284--294, 1990.

\bibitem{MeRo77}
C. D. Meyer and N. J. Rose.
\newblock The index and Drazin inverse of block triangular
matrices,
\newblock {\it SIAM J. Appl. Math.}, 33:1--7, 1977.


\bibitem{MoSc57}
J. L. Mott and H. Schneider.
\newblock Matrix norms applied to weakly ergodic Markov chains,
\newblock {\it Arch. Math.}, 8:331--333,  1957.

\bibitem{MoSc59}
J. L. Mott and H. Schneider.
\newblock Matrix algebras and groups relatively bounded in norm.
\newblock {\it Arch. Math.}, 10:1--6, 1959.

\bibitem{NeNe87}
S. Nelson and M. Neumann.
\newblock Generalization of the
Projection Method with Applications to SOR Method for Hermitian
Positive Semidefinite Linear Systems.
\newblock {\it Numer. Math.},  51:123--141, 1987.

\bibitem{Rhod97}
A. Rhodius,
On the maximum ergodicity coefficient, the Dobrushin coefficient
and products of stochastic matrices, {\it Lin. Alg. Appl.}
253:141--157, 1997.

\bibitem{Sene73}
E. Seneta,
On the historical development of the theory of finite
inhomogeneous
Markov chains,
\newblock{\it Proc. Cambridge Phil. Soc}, 74:507--513, 1973.


\bibitem{Sene81}
E. Seneta.
\newblock Non--negative Matrices and Markov Chains, Second
Edition.
\newblock {\it Springer Verlag}, New--York, 1981.

\bibitem{Zeng72}
C. Zenger, A comparison for some bounds for the nontrivial
eigenvalues
of stochastic matrices, {\it Numer. Math.}, 19:209--211, 1972.

\end{thebibliography}

\end{document}